\documentclass[10pt]{article}
\begin{document}
\title{ {\bf The smallest degree sum that
yields potentially $K_{r+1}-Z$-graphical Sequences}
\thanks{  Project Supported by NNSF of China(10271105), NSF of Fujian(Z0511034),
Fujian Provincial Training Foundation for "Bai-Quan-Wan Talents
Engineering" , Project of Fujian Education Department and Project of
Zhangzhou Teachers College.}}
\author{{Chunhui Lai}\\
{\small Department of Mathematics, Zhangzhou Teachers College,}
\\{\small Zhangzhou, Fujian 363000,
 P. R. of CHINA.}\\{\small e-mail: zjlaichu@public.zzptt.fj.cn }}
\date{}
\maketitle
\begin{center}
\begin{minipage}{4.3in}
\vskip 0.1in
\begin{center}{\bf Abstract}\end{center}
 { Let $K_{m}-H$ be the graph
obtained from $K_{m}$ by removing the edges set $E(H)$ of the graph
$H$ ($H$ is a subgraph of $K_{m}$). We use the symbol $Z_4$ to
denote $K_4-P_2.$  A sequence $S$ is potentially $K_{m}-H$-graphical
if it has a realization containing a $K_{m}-H$ as a subgraph. Let
$\sigma(K_{m}-H, n)$ denote the smallest degree sum such that every
$n$-term graphical sequence $S$ with $\sigma(S)\geq \sigma(K_{m}-H,
n)$ is potentially $K_{m}-H$-graphical.  In this paper, we determine
the values of $\sigma (K_{r+1}-Z, n)$ for
    $n\geq 5r+19,  r+1 \geq k \geq 5,$  $j \geq 5$ where $Z$ is a graph on $k$
    vertices and $j$ edges which
    contains a graph  $Z_4$  but
     not contains a cycle on $4$ vertices. We also determine the values of
      $\sigma (K_{r+1}-Z_4, n)$, $\sigma (K_{r+1}-(K_4-e), n)$,
      $\sigma (K_{r+1}-K_4, n)$ for
    $n\geq 5r+16, r\geq 4$.} \par
\par
\end{minipage}
\end{center}
 \par
 \section{Introduction}
\par

  The set of all non-increasing nonnegative integers sequence $\pi=$
  ($d_1 ,$ $d_2 ,$ $...,$ $d_n $) is denoted by $NS_n$.
  A sequence
$\pi\in NS_n$ is said to be graphic if it is the degree sequence of
a simple graph $G$ on $n$ vertices, and such a graph $G$ is called a
realization of $\pi$. The set of all graphic sequences in $NS_n$ is
denoted by $GS_n$. A graphical sequence $\pi$ is potentially
$H$-graphical if there is a realization of $\pi$ containing $H$ as a
subgraph, while $\pi$ is forcibly $H$-graphical if every realization
of $\pi$ contains $H$ as a subgraph. If $\pi$ has a realization in
which the $r+1$ vertices of  largest degree induce a clique, then
$\pi$ is said to be potentially $A_{r+1}$-graphic. Let
$\sigma(\pi)=d_1 +d_2 +... +d_n ,$ and $[x]$ denote the largest
integer less than or equal to $x$. If $G$ and $G_1$ are graphs, then
$G\cup G_1$ is the disjoint union of $G$ and $G_1$. If $G = G_1$, we
abbreviate $G\cup G_1$ as $2G$. We denote $G+H$ as the graph with
$V(G+H)=V(G)\bigcup V(H)$ and $E(G+H)=E(G)\bigcup E(H)\bigcup \{xy:
x\in V(G) , y \in V(H) \}. $ Let $K_k$, $C_k$, $T_k$, and $P_{k}$
denote a complete graph on $k$ vertices,  a cycle on $k$ vertices, a
tree on $k+1$ vertices, and a path on $k+1$ vertices, respectively.
Let $K_{m}-H$ be the graph obtained from $K_{m}$ by removing the
edges set $E(H)$ of the graph $H$ ($H$ is a subgraph of $K_{m}$). We
use the symbol $Z_4$ to denote $K_4-P_2.$  We use the symbol
$G[v_1,v_2,..., v_k]$ to denote the subgraph of $G$ induced by
vertex set $\{v_1,v_2,..., v_k \}$. We use the symbol $\epsilon (G)$
to denote the numbers of edges in graph $G$.
\par

Given a graph $H$, what is the maximum number of edges of a graph
with $n$ vertices not containing $H$ as a subgraph? This number is
denoted $ex(n,H)$, and is known as the Tur\'{a}n number. This
problem was proposed for $H = C_4$ by Erd\"os [2] in 1938 and in
general by Tur\'{a}n [19]. In terms of graphic sequences, the number
$2ex(n,H)+2$ is the minimum even integer $l$ such that every
$n$-term graphical sequence $\pi$ with $\sigma (\pi) \geq l $ is
forcibly $H$-graphical. Here we consider the following variant:
determine the minimum even integer $l$ such that every $n$-term
graphical sequence $\pi$ with $\sigma(\pi)\ge l$ is potentially
$H$-graphical. We denote this minimum $l$ by $\sigma(H, n)$.
Erd\"os,\ Jacobson and Lehel [4] showed that $\sigma(K_k, n)\ge
(k-2)(2n-k+1)+2$ and conjectured that equality holds. They proved
that if $\pi$ does not contain zero terms, this conjecture is true
for $k=3,\ n\ge 6$. The conjecture is confirmed in
[5],[14],[15],[16] and [17].
 \par
 Gould,\ Jacobson and
Lehel [5] also proved that  $\sigma(pK_2, n)=(p-1)(2n-2)+2$ for
$p\ge 2$; $\sigma(C_4, n)=2[{{3n-1}\over 2}]$ for $n\ge 4$. They
also pointed out that it would be nice to see where in the range for
$3n-2$  to $4n-4,$ the value  $\sigma (K_4-e, n)$ lies.  Luo [18]
characterized the potentially $C_{k}$ graphic sequence for
$k=3,4,5.$  Lai [7] determined  $\sigma (K_4-e, n)$ for $n\ge 4$.\
Yin,Li and Mao[21] determined $\sigma(K_{r+1}-e,n)$ for $r\geq 3,$
$r+1\leq n \leq 2r$ and $\sigma(K_5-e,n)$ for $n\geq5$. Yin and Li
[20] gave a good method (Yin-Li method) of determining the values
$\sigma(K_{r+1}-e,n)$ for $r\geq2$ and $n\geq3r^2-r-1$ (In fact, Yin
and Li[20] also determining the values $\sigma(K_{r+1}-ke,n)$ for
$r\geq2$ and $n\geq 3r^2-r-1$). After reading[20], using Yin-Li
method Yin [22] determined  $\sigma (K_{r+1}-K_{3}, n)$ for
    $n\geq 3r+5, r\geq 3$.  Lai [8] determined
    $\sigma (K_{5}-K_{3}, n),$ for  $n\geq 5$. Lai [9]
    gave a lower bound of $\sigma (K_{t+p}-K_{p}, n).$
      Lai [10,11] determined
    $\sigma (K_{5}-C_{4}, n),\sigma (K_{5}-P_{3}, n)$ and
    $\sigma (K_{5}-P_{4}, n),$ for  $n\geq 5$. Determining $\sigma(K_{r+1}-H,n)$, where $H$
    is a tree on 4 vertices is more useful than a cycle on 4
    vertices (for example, $C_4 \not\subset C_i$, but $P_3 \subset C_i$ for $i\geq 5$).
    So, after reading[20] and [22], using Yin-Li method Lai and Hu[12]
    determined  $\sigma (K_{r+1}-H, n)$ for
    $n\geq 4r+10, r\geq 3, r+1 \geq k \geq 4$ and $H$ be a graph on $k$
    vertices which
    containing a tree on $4$ vertices but
     not containing a cycle on $3$ vertices and $\sigma (K_{r+1}-P_2, n)$ for
    $n\geq 4r+8, r\geq 3$. Using Yin-Li method Lai and Sun[13] determined
$\sigma (K_{r+1}-(kP_2\bigcup tK_2), n)$ for
    $n\geq 4r+10, r+1 \geq 3k+2t,
    k+t \geq 2,k \geq 1, t \geq 0$ . To now, the problem of determining $\sigma
(K_{r+1}-H, n)$ for $H$ not containing a cycle on 3 vertices and
sufficiently large $n$ has been solved.
    In this paper, using Yin-Li method we prove the following two
theorems.\par

{\bf  Theorem 1.1.} If $r\geq 4$ and $n\geq 5r+16$, then
 $$ \sigma (K_{r+1}-K_{4}, n) = \sigma (K_{r+1}-(K_{4}-e), n) =$$
 $$\sigma (K_{r+1}-Z_{4}, n) =\left\{
    \begin{array}{ll}(r-1)(2n-r)-3(n-r)+1, \\ \mbox{ if $n-r$ is odd}\\
    (r-1)(2n-r)-3(n-r)+2,
     \\ \mbox{if $n-r$ is even} \end{array} \right. $$

\par
{\bf  Theorem 1.2.} If $n\geq 5r+19,  r+1 \geq k \geq 5,$ and $j
\geq 5$, then
  $$\sigma (K_{r+1}-Z, n) =\left\{
    \begin{array}{ll}(r-1)(2n-r)-3(n-r)-1, \\ \mbox{ if $n-r$ is odd}\\
    (r-1)(2n-r)-3(n-r)-2,
     \\ \mbox{if $n-r$ is even} \end{array} \right. $$
      where $Z$ is a graph on $k$
    vertices and $j$ edges which
    contains a graph $Z_4$ but
     not contains a cycle on $4$ vertices.
    \par

There are a number of graphs on $k$
    vertices  and $j$ edges which
    contains a graph $Z_4$  but
     not contains a cycle on $4$ vertices.
       \par

\section{Preparations}\par
  In order to prove our main result,we need the following notations
  and results.\par
  Let $\pi=(d_1,\cdots,d_n)\in NS_n,1\leq k\leq n$. Let \par
    $$ \pi_k^{\prime\prime}=\left\{
    \begin{array}{ll}(d_1-1,\cdots,d_{k-1}-1,d_{k+1}-1,
    \cdots,d_{d_k+1}-1,d_{d_k+2},\cdots,d_n), \\ \mbox{ if $d_k\geq k,$}\\
    (d_1-1,\cdots,d_{d_k}-1,d_{d_k+1},\cdots,d_{k-1},d_{k+1},\cdots,d_n),
     \\ \mbox{if $d_k < k.$} \end{array} \right. $$
  Denote
  $\pi_k^\prime=(d_1^\prime,d_2^\prime,\cdots,d_{n-1}^\prime)$,where
  $d_1^\prime\geq d_2^\prime\geq\cdots\geq d_{n-1}^\prime$ is a
  rearrangement of the $n-1$ terms of $\pi_k^{\prime\prime}$. Then
  $\pi_k^{\prime}$ is called the residual sequence obtained by
  laying off $d_k$ from $\pi$.\par
    {\bf Theorem 2.1[20]} Let $n\geq r+1$ and $\pi=(d_1,d_2,\cdots,d_n)\in
    GS_n$ with $d_{r+1}\geq r$. If $d_i\geq 2r-i$ for
    $i=1,2,\cdots,r-1$, then $\pi$ is potentially $A_{r+1}$-graphic.
    \par
    {\bf Theorem 2.2[20]} Let $n\geq 2r+2$ and $\pi=(d_1,d_2,\cdots,d_n)\in
    GS_n$ with $d_{r+1}\geq r$. If $d_{2r+2}\geq r-1$ , then $\pi$ is
    potentially $A_{r+1}$-graphic.
    \par
    {\bf Theorem 2.3[20]} Let $n\geq r+1$ and $\pi=(d_1,d_2,\cdots,d_n)\in
    GS_n$ with $d_{r+1}\geq r-1$. If $d_i\geq 2r-i$ for
    $i=1,2,\cdots,r-1$, then $\pi$ is potentially $K_{r+1}-e$-graphic.
    \par
    {\bf Theorem 2.4[20]} Let $n\geq 2r+2$ and $\pi=(d_1,d_2,\cdots,d_n)\in
    GS_n$ with $d_{r-1}\geq r$. If $d_{2r+2}\geq r-1$ , then $\pi$ is
    potentially $K_{r+1}-e$
    -graphic.
    \par
    {\bf Theorem 2.5[6]} Let $\pi=(d_1,\cdots,d_n)\in NS_n$ and $1\leq k\leq
    n$. Then $\pi\in GS_n$ if and only if  $\pi_k^\prime\in
    GS_{n-1}$.
    \par
     {\bf Theorem 2.6[3]} Let $\pi=(d_1,\cdots,d_n)\in NS_n$
     with even $\sigma(\pi)$. Then $\pi\in GS_n$ if and only if
     for any $t$,$1\leq t\leq n-1$,
     $$\sum_{i=1}^t d_i\leq t(t-1)+\sum_{j=t+1}^n
     min \{t,d_j \}.$$
     \par
     {\bf Theorem 2.7[5]} If $\pi=(d_1,d_2,\cdots,d_n)$ is a graphic
     sequence with a realization $G$ containing $H$ as a subgraph, then
     there exists a realization $G^\prime$ of $\pi$ containing H as a
     subgraph so that the vertices of $H$ have the largest degrees of
     $\pi$.
     \par
     {\bf Theorem 2.8[9]} If $n \geq p+t,$
     then $\sigma (K_{p+t}-K_{p}, n)\geq 2[((p+2t-3)n+p+2t+1-pt-t^2)/2]$.
\par
     \par
     {\bf Lemma 2.1 [22]} If $\pi=(d_1,d_2,\cdots,d_n)\in NS_n$ is potentially
     $K_{r+1}-e$-graphic, then there is a realization $G$ of $\pi$
     containing $K_{r+1}-e$ with the $r+1$ vertices
     $v_1,\cdots,v_{r+1}$ such that $d_G(v_i)=d_i$ for
     $i=1,2,\cdots,r+1$ and $e=v_rv_{r+1}$.
     \par
     {\bf Lemma 2.2 [12]}  Let $n\geq 2r+2$ and $\pi=(d_1,d_2,\cdots,d_n)\in GS_n$
with  $d_{r-2}\geq r$. If $d_{2r+2}\geq r-1$, then $\pi$ is
potentially $K_{r+1}-P_2$-graphic.
\par
     {\bf Lemma 2.3} Let $\pi=(d_1,\cdots,d_n)\in GS_n$
     and $G$ be a realization of $\pi$. If $\epsilon (G[v_1,v_2,..., v_{r+1}]) \leq \epsilon
     (K_{r+1})- 1$, then there is a realization $H$ of $\pi$
      such that $d_H(v_i)=d_i$ for
     $i=1,2,\cdots,r+1$ and $v_rv_{r+1}\not\in E(H)$.
     \par
     The proof is similar to the proof of Lemma 2.1.
     \par

\section{ Proof of Main results.} \par

{\bf  Lemma 3.1.} Let $n\geq 2r$ and $\pi=(d_1,d_2,\cdots,d_n)\in
    GS_n$ with $d_{r-1}\geq r$, $d_{r+1}\geq r-1$. If $d_i\geq 2r-i$ for
    $i=1,2,\cdots,r-2$, then $\pi$ is potentially $K_{r+1}-e$-graphic.
    \par
  {\bf Proof.}   We consider the following two cases.
  \par
  Case 1: $d_{r+1}\geq r.$ \par
  If $d_{r-1}\geq r+1.$ \par
  Then $\pi$ is potentially
  $K_{r+1}-e$-graphic by Theorem 2.3.
  \par
  If $d_{r-1}= r,$ then $d_{r-1}=d_{r}=d_{r+1}=r$ \par
  Suppose $\pi$ is not potentially $K_{r+1}-e$-graphic.
  Let $H$ be a realization of $\pi$, then $\epsilon
  (H[v_1,v_2,..., v_{r+1}]) \leq \epsilon
     (K_{r+1})- 2.$
  Let $S =(d_1,d_2,\cdots,d_{r-2},
  d_{r-1},$ $d_{r}+1,d_{r+1}+1,$ $ \cdots, d_n)$, then by Theorem 2.1,
   $S$ is potentially   $A_{r+1}$-graphic (Denote
  $S^\prime=(d_1^\prime,d_2^\prime,$ $\cdots,d_{n}^\prime)$,where
  $d_1^\prime\geq d_2^\prime\geq\cdots\geq d_{n}^\prime$ is a
  rearrangement of the $n$ terms of $S$. Therefore $S^\prime \in GS_n$
  by  Lemma 2.3. Then  $S^\prime$
  satisfies the conditions of Theorem 2.1). Therefore, there is a
   realization $G$ of $S$ with $v_{1},v_{2},\cdots, v_{r+1}$
   $(d(v_i)=d_i, i=1,2,\cdots,r-1,$
   $d(v_{r})=d_{r}+1,d(v_{r+1})=d_{r+1}+1 ),$  the $r+1$
   vertices of highest degree containing a $K_{r+1}$. Hence,
   $G-v_{r+1}v_{r}$ is a  realization of $\pi$.
   Thus,  $\pi$ is potentially $K_{r+1}-e$-graphic,
   which is a contradiction.
  \par
  Case 2: $d_{r+1} = r-1$,  then the residual sequence $\pi_{r+1}^\prime=
  (d_1^\prime,\cdots,d_{n-1}^\prime)$
  obtained by laying off $d_{r+1}=r-1$ from $\pi$ satisfies:
  $d_1^\prime\geq2(r-1)-1$,$ \cdots,$
  $d_{(r-1)-1}^\prime=d_{r-2}^\prime\geq 2(r-1)-(r-2), d_{(r-1)+1}^\prime=d_r^\prime
   \geq r-1$. By Theorem 2.1, $\pi_{r+1}^\prime$ is potentially
  $A_{(r-1)+1}$-graphic. Therefore, $\pi$ is potentially
  $K_{r+1}-e$-graphic by $ \{d_1-1,\cdots,d_{r-1}-1 \}
  \subseteq  \{d_1^\prime,\cdots,d_r^\prime \}$
  and Theorem 2.7.
  \par
  {\bf  Lemma 3.2.} Let $n\geq 2r$ and
$\pi=(d_1,d_2,\cdots,d_n)\in
    GS_n$ with $d_{r-2}\geq r+1$, $d_{r+1} \geq r$,$d_{r}-1 \geq d_{d_{r+1}+2}$. If $d_i\geq 2r-i$ for
    $i=1,2,\cdots,r-3$, then $\pi$ is potentially $A_{r+1}$-graphic.
    \par
  {\bf Proof.}
    The residual sequence $\pi_{r+1}^\prime=
  (d_1^\prime,\cdots,d_{n-1}^\prime)$
  obtained by laying off $d_{r+1}$ from $\pi$ satisfies:
  $d_1^\prime\geq2(r-1)-1$,$ \cdots,$
  $d_{(r-1)-2}^\prime=d_{r-3}^\prime\geq 2(r-1)-(r-3),$ $d_{(r-1)-1}^\prime=d_{r-2}^\prime\geq 2(r-1)-(r-2),$
   $d_{(r-1)+1}^\prime=d_r^\prime
   \geq r-1$. By Theorem 2.1, $\pi_{r+1}^\prime$ is potentially
  $A_{(r-1)+1}$-graphic. Therefore, $\pi$ is potentially
  $A_{r+1}$-graphic by $ \{d_1-1,\cdots,d_{r}-1 \}
  =  \{d_1^\prime,\cdots,d_r^\prime \}$
  and Theorem 2.7.
  \par
   {\bf Lemma 3.3 } Let $n\geq 2r+2, r\geq  4$ and $\pi=(d_1,d_2,\cdots,d_n)\in GS_n$
with $d_{r-2}\geq r-1$ and $d_{r+1}\geq r-2$,
$$\sigma (\pi) \geq \left\{
    \begin{array}{ll}(r-1)(2n-r)-3(n-r)-1, \\ \mbox{ if $n-r$ is odd}\\
    (r-1)(2n-r)-3(n-r)-2,
     \\ \mbox{if $n-r$ is even} \end{array} \right. $$
      If $d_i\geq 2r-i$ for
$i=1,2,\cdots,r-3,$ then $\pi$ is potentially $K_{r+1}-Z_4$-graphic.

    \par

  {\bf Proof.} We consider the following two cases.
  \par
  Case 1: $d_{r+1}\geq r-1.$
  \par
  Subcase 1.1:  $d_{r-1}\geq r+1$. \par
  If $d_{r-2}\geq r+2$, then $\pi$ is potentially
  $K_{r+1}-e$-graphic by Theorem 2.3.
  Hence,$\pi$ is potentially $K_{r+1}-Z_4$-graphic.
  \par
  If $d_{r-2}=r+1,$ then $d_{r-3}-1\geq d_{r-2}$.
   The residual sequence $\pi_{r+1}^\prime=
  (d_1^\prime,\cdots,d_{n-1}^\prime)$
  obtained by laying off $d_{r+1}$ from $\pi$ satisfies:
  $d_1^\prime\geq2(r-1)-1$,$ \cdots,$
  $d_{(r-1)-2}^\prime=d_{r-3}^\prime\geq 2(r-1)-(r-3),$
  $d_{(r-1)-1}^\prime=d_{r-2}^\prime\geq r-1,$
   $d_{(r-1)+1}^\prime=d_r^\prime
   \geq (r-1)-1$. By Lemma 3.1, $\pi_{r+1}^\prime$ is potentially
  $K_{(r-1)+1}-e$-graphic. Therefore, $\pi$ is potentially
  $K_{r+1}-Z_4$-graphic by $ \{d_1-1,\cdots,d_{r-3}-1 \}
  \subseteq  \{d_1^\prime,\cdots,d_r^\prime \}$
  and Lemma 2.1.

    \par

    Subcase 1.2:  $d_{r-1}\leq r$. then $d_{r-3}-1\geq d_{r-1}$.
   The residual sequence $\pi_{r+1}^\prime=
  (d_1^\prime,\cdots,d_{n-1}^\prime)$
  obtained by laying off $d_{r+1}$ from $\pi$ satisfies:
  $d_1^\prime\geq2(r-1)-1$,$ \cdots,$
  $d_{(r-1)-2}^\prime=d_{r-3}^\prime\geq 2(r-1)-(r-3),$
  $d_{(r-1)-1}^\prime=d_{r-2}^\prime\geq r-1,$
   $d_{(r-1)+1}^\prime=d_r^\prime
   \geq (r-1)-1$. By Lemma 3.1, $\pi_{r+1}^\prime$ is potentially
  $K_{(r-1)+1}-e$-graphic. Therefore, $\pi$ is potentially
  $K_{r+1}-Z_4$-graphic by $ \{d_1-1,\cdots,d_{r-3}-1 \}
  \subseteq  \{d_1^\prime,\cdots,d_r^\prime \}$
  and Lemma 2.1.
    \par

  Case 2:  $d_{r+1}= r-2$.
\par
If $d_{r-1}<d_{ r-2}.$
 \par
 If $d_{r-2}\geq r,$
 then the residual sequence $\pi_{r+1}^\prime=
  (d_1^\prime,\cdots,d_{n-1}^\prime)$
  obtained by laying off $d_{r+1}=r-2$ from $\pi$ satisfies: (1)
$d_i^\prime=d_i-1$ for
  $i=1,2,\cdots,r-2$,(2) $d_1^\prime=d_1-1\geq 2(r-1)-1,
  \cdots,d_{(r-1)-2}^\prime=d_{r-3}^\prime \geq d_{r-3}-1
  \geq 2(r-1)-[(r-1)-2]$, $d_{(r-1)-1}^\prime=d_{r-2}^\prime
  \geq r-1,$ and $d_{(r-1)+1}^\prime=d_r^\prime=d_r\geq r-2$.
  By Lemma 3.1, $\pi_{r+1}^\prime$ is potentially
  $K_{(r-1)+1}-e$-graphic. Therefore, $\pi$ is potentially
  $K_{r+1}-Z_4$-graphic by $\{d_1-1,\cdots,d_{r-2}-1,d_{r-1},d_r \}=
  \{d_1^\prime,\cdots,d_r^\prime \}$
  and Lemma 2.1.
  \par
 If $d_{r-2}= r-1,$ then $d_{r-1}=d_{r}=r-2$ and
 $$ \begin{array}{rcl}
   \sigma(\pi)&\leq &(r-3)(n-1)+r-1+(r-2)(n-r+2)\\
   &=&   (r-1)(n-1)-2(n-1)+(r-1)(n-r+3)-(n-r+2)\\
   &=& (r-1)(2n-r)-3(n-r)-2\\
   \end{array} $$
   Hence, $\pi=( (n-1)^{r-3}, (r-1)^1,(r-2)^{n-r+2})$ and $n-r$ is even.
   Clearly,  $\pi$ is potentially $K_{r+1}-Z_4$-graphic.
    \par
  \par
  If $d_{r-1}=d_{ r-2}$ and $d_{ r-3}\geq d_{ r}$, then $\pi_{r+1}^\prime$ satisfies:
   $d_1^\prime \geq d_1-1 \geq 2(r-1)-1,\cdots,d_{(r-1)-2}^\prime=
   d_{r-3}^\prime \geq d_{r-3}-1
  \geq 2(r-1)-[(r-1)-2]$, $d_{(r-1)-1}^\prime=
   d_{r-2}^\prime
  \geq r-1$ and $d_{(r-1)+1}^\prime=d_r^\prime \geq r-2$.
  By Lemma 3.1, $\pi_{r+1}^\prime$ is potentially
  $K_{(r-1)+1}-e$-graphic. Therefore, $\pi$ is potentially
  $K_{r+1}-Z_4$-graphic by $\{d_{r-1},d_r,d_1-1\cdots,d_{r-2}-1 \}=
  \{d_1^\prime,
  \cdots,d_r^\prime \}$
  and Lemma 2.1.
  \par
  If $d_{r-1}=d_{ r-2}$ and $d_{ r-3}= d_{ r}$, then $d_{ r-3}=d_{r-2}=d_{ r-1}= d_{
  r}\geq r+3$. Let $H$ be a realization of $\pi$. Since $d_{r+1}=
  r-2$,  then there is $i,j\leq r$ such that $v_{r+1}v_{i}, v_{r+1}v_{j} \not\in E(H).$
 Let $S =(d_1,d_2,\cdots,d_{i}+1,\cdots, d_{j}+1,\cdots,$
 $d_{r},d_{r+1}+2,$ $ \cdots, d_n)$, then by Theorem 2.1,
   $S$ is potentially   $A_{r+1}$-graphic (Denote
  $S^\prime=(d_1^\prime,d_2^\prime,$ $\cdots,d_{n}^\prime)$,where
  $d_1^\prime\geq d_2^\prime\geq\cdots\geq d_{n}^\prime$ is a
  rearrangement of the $n$ terms of $S$. Therefore $S^\prime \in GS_n$.
  Then  $S^\prime$
  satisfies the conditions of Theorem 2.1). Therefore, there is a
   realization $G$ of $S$ with $v_{1},v_{2},\cdots, v_{r+1}$
   $(d(v_t)=d_t, t\neq i, j, r+1,$
   $d(v_{i})=d_{i}+1,d(v_{j})=d_{j}+1,d(v_{r+1})=d_{r+1}+2 ),$  the $r+1$
   vertices of highest degree containing a $K_{r+1}$. Hence,
   $G-\{v_{r+1}v_{i},v_{r+1}v_{j}\}$ is a  realization of $\pi$.
   Thus,  $\pi$ is potentially $K_{r+1}-Z_4$-graphic.
\par
{\bf Lemma 3.4 } Let $n\geq 2r+2$ and $\pi=(d_1,d_2,\cdots,d_n)\in
GS_n$ with  $d_{r-t}\geq r$. If $d_{2r+2}\geq r-1$, then $\pi$ is
potentially $K_{r+1}-K_{1,t}$-graphic.
\par
{\bf  Proof. }
  We consider the following two cases.
  \par
  Case 1: If $d_{r-1}\geq r$. Then $\pi$ is potentially
  $K_{r+1}-e$-graphic by Theorem 2.4. Hence, $\pi$ is
potentially $K_{r+1}-K_{1,t}$-graphic.
  \par
  Case 2: $d_{r-1}\leq r-1$, that is, $d_{r-1}= r-1$, then
  $d_{r-1}=d_r=d_{r+1}=\cdots=d_{2r+2}=r-1$ and $\pi_{r+1}^\prime$ satisfies:
   $d_{(r-1)+1}^\prime=d_r^\prime\geq r-1$ and
   $d_{2(r-1)+2}^\prime=d_{2r}^\prime\geq (r-1)-1$.
  By Theorem 2.2, $\pi_{r+1}^\prime$ is potentially
  $A_r$-graphic. Therefore, $\pi$ is potentially
  $K_{r+1}-K_{1,t}$-graphic by $\{d_1-1,\cdots,d_{r-t}-1 \}
  \subseteq \{d_1^\prime,\cdots,d_r^\prime \}$
  and Theorem 2.7.
  \par
  {\bf Lemma 3.5 } Let $n\geq 2r+2$ and
$\pi=(d_1,d_2,\cdots,d_n)\in GS_n$ with  $d_{r-4}\geq r$,
$$\sigma (\pi) \geq \left\{
    \begin{array}{ll}(r-1)(2n-r)-3(n-r)-1, \\ \mbox{ if $n-r$ is odd}\\
    (r-1)(2n-r)-3(n-r)-2,
     \\ \mbox{if $n-r$ is even} \end{array} \right. $$
      If
$d_{2r+2}\geq r-1$, then $\pi$ is potentially
 $K_{r+1}-(P_{2}\bigcup K_{2})$-graphic.
\par
{\bf  Proof. }
  We consider the following two cases.
  \par
  Case 1: If $d_{r-2}\geq r$. Then $\pi$ is potentially
  $K_{r+1}-P_{2}$-graphic by Lemma 2.2. Hence, $\pi$ is
potentially $K_{r+1}-(P_{2}\bigcup K_{2})$-graphic.
  \par
  Case 2: $d_{r-2}= r-1$. \par
Subcase 2.1: $d_{r-3}\geq r$, then $d_{r-3}\geq
d_{r}+1=d_{r+1}+1=r>r-1=d_{r-2}=d_{r-1}$. Suppose $\pi$ is not
potentially $K_{r+1}-(P_{2}\bigcup K_{2})$-graphic.
  Let $H$ be a realization of $\pi$, then $\epsilon
  (H[v_1,v_2,..., v_{r+1}]) \leq \epsilon
     (K_{r+1})- 3.$
     Let $S =(d_1,d_2,\cdots,d_{r-2},
  d_{r-1},d_{r}+1,d_{r+1}+1,$ $ \cdots, d_n)$, then by Theorem 2.4,
   $S$ is potentially   $K_{r+1}-e$-graphic (Denote
  $S^\prime=(d_1^\prime,d_2^\prime,$ $\cdots,d_{n}^\prime)$,where
  $d_1^\prime\geq d_2^\prime\geq\cdots\geq d_{n}^\prime$ is a
  rearrangement of the $n$ terms of $S$. Therefore $S^\prime \in GS_n$
  by  Lemma 2.3. Then  $S^\prime$
  satisfies the conditions of Theorem 2.4). Therefore, there is a
   realization $G$ of $S$ with $v_{1},v_{2},\cdots, v_{r+1}$
   $(d(v_i)=d_i, i=1,2,\cdots,r-1,$
   $d(v_{r})=d_{r}+1,d(v_{r+1})=d_{r+1}+1 ),$  the $r+1$
   vertices of highest degree containing a $K_{r+1}-e$
   and $e=v_{r-1}v_{r-2}$ by Lemma 2.1. Hence,
   $G-v_{r+1}v_{r}$ is a  realization of $\pi$.
   Thus,  $\pi$ is potentially $K_{r+1}-(P_{2}\bigcup
   K_{2})$-graphic,
   which is a contradiction.
\par
Subcase 2.2: $d_{r-3}= r-1$, then
  $$ \begin{array}{rcl}
   \sigma(\pi)&\leq &(r-4)(n-1)+(r-1)(n-r+4)\\
   &=&   (r-1)(n-1)-3(n-1)+(r-1)(n-r+1)+3(r-1)\\
   &=& (r-1)(2n-r)-3(n-r)\\
   \end{array} $$
   Since,
   $$\sigma (\pi) \geq \left\{
    \begin{array}{ll}(r-1)(2n-r)-3(n-r)-1, \\ \mbox{ if $n-r$ is odd}\\
    (r-1)(2n-r)-3(n-r)-2,
     \\ \mbox{if $n-r$ is even} \end{array} \right. $$
   Hence,  $\pi$ is one of the following:
   $( (n-1)^{r-5}, (n-2)^1,(r-1)^{n-r+4})$,
   $( (n-1)^{r-4}, (r-1)^{n-r+3}, (r-2)^1)$, for $n-r$ is odd,
     $\pi$ is one of the following:
    $( (n-1)^{r-4}, (r-1)^{n-r+4} )$,
    $( (n-1)^{r-6}, (n-2)^2,(r-1)^{n-r+4})$,
     $( (n-1)^{r-5}, (n-3)^1,(r-1)^{n-r+4})$,
      $( (n-1)^{r-5}, (n-2)^1,(r-1)^{n-r+3}, (r-2)^1)$,
      $( (n-1)^{r-4}, (r-1)^{n-r+3}, (r-3)^1)$,
      $( (n-1)^{r-4}, (r-1)^{n-r+2}, (r-2)^2)$,
      for  $n-r$ is even.
   Clearly,  $\pi$ is potentially $K_{r+1}-(P_{2}\bigcup
   K_{2})$-graphic.

\par
   {\bf Lemma 3.6.} If $r\geq 4$ and $n\geq r+1$, then
   $$\sigma(K_{r+1}-Z_4,n)\geq
    \sigma(K_{r+1}-K_4,n).$$
    and
          $$ \sigma (K_{r+1}-K_{4}, n) \geq\left\{
    \begin{array}{ll}(r-1)(2n-r)-3(n-r)+1, \\ \mbox{ if $n-r$ is odd}\\
    (r-1)(2n-r)-3(n-r)+2,
     \\ \mbox{if $n-r$ is even} \end{array} \right. $$
          \par
    {\bf Proof.} Obviously, for $r\geq 4$ and $n\geq r+1$,
      $\sigma(K_{r+1}-Z_4,n)\geq
    \sigma(K_{r+1}-K_4,n).$ By Theorem 2.8, for $r\geq 4$ and $n\geq r+1$,
    $\sigma (K_{r+1}-K_{4}, n)= \sigma (K_{4+(r-3)}-K_{4}, n)$
    $\geq 2[((4+2(r-3)-3)n+4+2(r-3)+1-4(r-3)-(r-3)^2)/2]$.
     Hence,
    $$ \sigma (K_{r+1}-K_{4}, n) \geq\left\{
    \begin{array}{ll}(r-1)(2n-r)-3(n-r)+1, \\ \mbox{ if $n-r$ is odd}\\
    (r-1)(2n-r)-3(n-r)+2,
     \\ \mbox{if $n-r$ is even} \end{array} \right. $$
    \par

    {\bf  Lemma 3.7.} If $n\geq r+1, r+1 \geq k \geq 4,$
     then
     $$ \sigma (K_{r+1}-H, n) \geq\left\{
    \begin{array}{ll}(r-1)(2n-r)-3(n-r)-1, \\ \mbox{ if $n-r$ is odd}\\
    (r-1)(2n-r)-3(n-r)-2,
     \\ \mbox{if $n-r$ is even} \end{array} \right. $$
      where $H$ is a graph on $k$
    vertices which not contains a cycle on $4$ vertices.
\par
{\bf Proof.}   Let
 $$ G=\left\{
    \begin{array}{ll}K_{r-3}+(\frac{n-r+1}{2}+1)K_{2}, \\ \mbox{ if $n-r$ is odd}\\
    K_{r-3}+(\frac{n-r+2}{2}K_{2}\bigcup K_{1}),
     \\ \mbox{if $n-r$ is even} \end{array} \right. $$
Then $G$ is a unique realization of
$$ \pi=\left\{
    \begin{array}{ll}((n-1)^{r-3}, (r-2)^{n-r+3}), \\ \mbox{ if $n-r$ is odd}\\
    ((n-1)^{r-3}, (r-2)^{n-r+2}, (r-3)^1),
     \\ \mbox{if $n-r$ is even} \end{array} \right. $$
and $G$ clearly does not contain $K_{r+1}-H$, where the symbol $x^y$
means $x$ repeats $y$ times in the sequence. Thus  $\sigma
(K_{r+1}-H, n) \geq \sigma(\pi)+2$. Therefore,
$$ \sigma (K_{r+1}-H, n) \geq\left\{
    \begin{array}{ll}(r-1)(2n-r)-3(n-r)-1, \\ \mbox{ if $n-r$ is odd}\\
    (r-1)(2n-r)-3(n-r)-2,
     \\ \mbox{if $n-r$ is even} \end{array} \right. $$
\par

 {\bf The Proof of Theorem 1.1 }
      According to Lemma 3.6 and $ \sigma (K_{r+1}-K_{4}, n)
       \leq \sigma (K_{r+1}-(K_{4}-e), n) \leq
 \sigma (K_{r+1}-Z_{4}, n),$ it is enough to verify that for
        $n\geq 5r+16$,
       $$ \sigma (K_{r+1}-Z_{4}, n) \leq \left\{
    \begin{array}{ll}(r-1)(2n-r)-3(n-r)+1, \\ \mbox{ if $n-r$ is odd}\\
    (r-1)(2n-r)-3(n-r)+2,
     \\ \mbox{if $n-r$ is even} \end{array} \right. $$
    We now prove that if $n\geq 5r+16$ and $\pi=(d_1,d_2,\cdots,d_n)\in GS_n$
with
$$ \sigma (\pi) \geq \left\{
    \begin{array}{ll}(r-1)(2n-r)-3(n-r)+1, \\ \mbox{ if $n-r$ is odd}\\
    (r-1)(2n-r)-3(n-r)+2,
     \\ \mbox{if $n-r$ is even} \end{array} \right. $$
then $\pi$ is potentially
  $K_{r+1}-Z_4$-graphic.
  \par
   If $d_{r-3}\leq r-1$, then
 $$ \begin{array}{rcl}
   \sigma(\pi)&\leq &(r-4)(n-1)+(r-1)(n-r+4)\\
   &=&   (r-1)(n-1)-3(n-1)+(r-1)(n-r+4)\\
   &=& (r-1)(2n-r)-3(n-r)\\
   & < & (r-1)(2n-r)-3(n-r)+1,
   \end{array} $$
 which is
 a contradiction.  Thus, $d_{r-3}\geq r.$
 \par
 If $d_{r-2}\leq r-2$, then
 $$ \begin{array}{rcl}
   \sigma(\pi)&\leq &(r-3)(n-1)+(r-2)(n-r+3)\\
   &=&   (r-1)(n-1)-2(n-1)+(r-1)(n-r+3)-(n-r+3)\\
   &=& (r-1)(2n-r)-3(n-r)-3\\
   & < & (r-1)(2n-r)-3(n-r)+1,
   \end{array} $$
 which is
 a contradiction.  Thus, $d_{r-2}\geq r-1.$

\par
    If $d_{r+1}\leq r-3,$
    then
    $$\begin{array}{rcl}
    \sigma(\pi)&=&\sum_{i=1}^{r}d_i+\sum_{i=r+1}^n d_i \\
    & \leq &(r-1)r+
    \sum_{i=r+1}^{n}min \{r,d_i \}+\sum_{i=r+1}^n d_i \\
    &=&(r-1)r+2\sum_{i=r+1}^n
    d_i \\
    & \leq &(r-1)r+2(n-r)(r-3) \\
    &=&(r-1)(2n-r)-4(n-r) \\
    &<&(r-1)(2n-r)-3(n-r)+1, \end{array} $$
     which is a
    contradiction. Thus, $d_{r+1}\geq r-2$.
    \par
      If $d_i\geq 2r-i$ for $i=1,2,\cdots,r-3$ or $d_{2r+2}\geq r-1$, then
       $\pi$
      is potentially
  $K_{r+1}-Z_4$-graphic  by Lemma 3.3 or Lemma 3.4.  If $d_{2r+2}\leq
  r-2$ and there exists an integer $i$, $1\leq i\leq r-3$ such that
  $d_i\leq 2r-i-1$, then
  $$\begin{array}{rcl}
  \sigma(\pi) &\leq &(i-1)(n-1)+(2r+1-i+1)(2r-i-1)\\
  &&+
  (r-2)(n+1-2r-2)\\
  &=& i^2+i(n-4r-2)-(n-1)\\
  &&+(2r-1)(2r+2)+(r-2)(n-2r-1).
  \end{array} $$
  Since $n\geq 5r+16$, it is  easy to see that
  $i^2+i(n-4r-2)$, consider as a function of $i$, attains its maximum
  value when $i=r-3$. Therefore,
  $$\begin{array}{rcl}
  \sigma(\pi) &\leq &
  (r-3)^2+(n-4r-2)(r-3)-(n-1)\\
  &&+(2r-1)(2r+2)+(r-2)(n-2r-1)\\
  &=&(r-1)(2n-r)-3(n-r)-n+5r+16 \\
  &<&\sigma(\pi),
  \end{array} $$
  which is a
    contradiction.
    \par
    Thus,
     $$ \sigma (K_{r+1}-Z_{4}, n) \leq \left\{
    \begin{array}{ll}(r-1)(2n-r)-3(n-r)+1, \\ \mbox{ if $n-r$ is odd}\\
    (r-1)(2n-r)-3(n-r)+2,
     \\ \mbox{if $n-r$ is even} \end{array} \right. $$
      for  $n\geq 5r+16$.
\par
{\bf The Proof of Theorem 1.2 }
       According to Lemma 3.7, it is enough to verify that for
        $n\geq 5r+19$,
       $$ \sigma (K_{r+1}-Z, n) \leq \left\{
    \begin{array}{ll}(r-1)(2n-r)-3(n-r)-1, \\ \mbox{ if $n-r$ is odd}\\
    (r-1)(2n-r)-3(n-r)-2,
     \\ \mbox{if $n-r$ is even} \end{array} \right. $$
    We now prove that if $n\geq 5r+19$ and $\pi=(d_1,d_2,\cdots,d_n)\in GS_n$
with
$$ \sigma (\pi) \geq \left\{
    \begin{array}{ll}(r-1)(2n-r)-3(n-r)-1, \\ \mbox{ if $n-r$ is odd}\\
    (r-1)(2n-r)-3(n-r)-2,
     \\ \mbox{if $n-r$ is even} \end{array} \right. $$
then $\pi$ is potentially
  $K_{r+1}-Z$-graphic.
  \par
   If $d_{r-4}\leq r-1$, then
 $$ \begin{array}{rcl}
   \sigma(\pi)&\leq &(r-5)(n-1)+(r-1)(n-r+5)\\
   &=&   (r-1)(n-1)-4(n-1)+(r-1)(n-r+5)\\
   &=& (r-1)(2n-r)-4(n-r)\\
   & < & (r-1)(2n-r)-3(n-r)-2,
   \end{array} $$
 which is
 a contradiction.  Thus, $d_{r-4}\geq r.$
 \par
 If $d_{r-2}\leq r-2$, then
 $$ \begin{array}{rcl}
   \sigma(\pi)&\leq &(r-3)(n-1)+(r-2)(n-r+3)\\
   &=&   (r-1)(n-1)-2(n-1)+(r-1)(n-r+3)-(n-r+3)\\
   &=& (r-1)(2n-r)-3(n-r)-3\\
   & < & (r-1)(2n-r)-3(n-r)-2,
   \end{array} $$
 which is
 a contradiction.  Thus, $d_{r-2}\geq r-1.$
    \par
     If $d_{r+1}\leq r-3,$
    then
    $$\begin{array}{rcl}
    \sigma(\pi)&=&\sum_{i=1}^{r}d_i+\sum_{i=r+1}^n d_i \\
    & \leq &(r-1)r+
    \sum_{i=r+1}^{n}min \{r,d_i \}+\sum_{i=r+1}^n d_i \\
    &=&(r-1)r+2\sum_{i=r+1}^n
    d_i \\
        & \leq &(r-1)r+2(n-r)(r-3) \\
    &=&(r-1)(2n-r)-4(n-r) \\
    &<&(r-1)(2n-r)-3(n-r)-2, \end{array} $$
     which is a
    contradiction. Thus, $d_{r+1}\geq r-2$.
    \par
      If $d_i\geq 2r-i$ for $i=1,2,\cdots,r-3$ or $d_{2r+2}\geq r-1$, then
       $\pi$
      is potentially
  $K_{r+1}-Z$-graphic  by Lemma 3.3 or Lemma 3.5 .  If $d_{2r+2}\leq
  r-2$ and there exists an integer $i$, $1\leq i\leq r-3$ such that
  $d_i\leq 2r-i-1$, then
  $$\begin{array}{rcl}
  \sigma(\pi) &\leq &(i-1)(n-1)+(2r+1-i+1)(2r-i-1)\\
  &&+
  (r-2)(n+1-2r-2)\\
  &=& i^2+i(n-4r-2)-(n-1)\\
  &&+(2r-1)(2r+2)+(r-2)(n-2r-1).
  \end{array} $$
  Since $n\geq 5r+19$, it is  easy to see that
  $i^2+i(n-4r-2)$, consider as a function of $i$, attains its maximum
  value when $i=r-3$. Therefore,
  $$\begin{array}{rcl}
  \sigma(\pi) &\leq &
  (r-3)^2+(n-4r-2)(r-3)-(n-1)\\
  &&+(2r-1)(2r+2)+(r-2)(n-2r-1)\\
  &=&(r-1)(2n-r)-3(n-r)-n+5r+16 \\
  &<&\sigma(\pi),
  \end{array} $$
  which is a
    contradiction.
    \par
    Thus,
     $$ \sigma (K_{r+1}-Z, n) \leq \left\{
    \begin{array}{ll}(r-1)(2n-r)-3(n-r)-1, \\ \mbox{ if $n-r$ is odd}\\
    (r-1)(2n-r)-3(n-r)-2,
     \\ \mbox{if $n-r$ is even} \end{array} \right. $$
      for  $n\geq 5r+19$.
\par

 \section*{Acknowledgment}
  The authors thanks the referees for many helpful comments.

\par

\end{document}